# ON LAWS OF LARGE NUMBERS FOR RANDOM WALKS

By Anders Karlsson[1] and François Ledrappier[2]

*Royal Institute of Technology and University of Notre Dame*

We prove a general noncommutative law of large numbers. This applies in particular to random walks on any locally finite homogeneous graph, as well as to Brownian motion on Riemannian manifolds which admit a compact quotient. It also generalizes Oseledec's multiplicative ergodic theorem. In addition, we show that $\varepsilon$-shadows of any ballistic random walk with finite moment on any group eventually intersect. Some related results concerning Coxeter groups and mapping class groups are recorded in the last section.

**1. Introduction.** The strong law of large numbers, certainly a fundamental result in probablility theory, asserts that for a sequence of i.i.d. random variables $X_i$ taking values in $\mathbf{R}$,

$$\frac{1}{n}(X_1 + X_2 + \cdots + X_n) \to E(X_1)$$

almost surely, provided $E|X_1| < \infty$. One might wonder, as Bellman [5], Kesten [22] and Furstenberg [11] did in the 1950's and 60's, whether there exist generalizations of this law when the random variables instead take values in a more general, noncommutative group $G$. In other words, what can be said about the behavior of

$$Z_n := g_1 g_2 \cdots g_n$$

as $n \to \infty$ where $g_i$ are i.i.d. in $G$ (cf. the Introduction of [11])? Note that it is not even clear how to formulate a statement that would generalize the law of large numbers, let alone how to prove it. The aforementioned three

Received April 2005; revised December 2005.
[1]Supported by Swedish Research Council (VR) Grant 2002-4771.
[2]Supported in part by NSF Grant DMS-05-00630.

authors studied the case of free groups and groups of matrices. For example, Furstenberg and Kesten proved in [10] that for matrices, the limit

$$\lim_{n\to\infty} \frac{1}{n} \log \|Z_n\|,$$

where $\|\cdot\|$ denotes a matrix norm, exists almost surely. Later on in that decade, this was generalized in two very important ways: first, Kingman's subadditive ergodic theorem [23], which showed that the reason for the above convergence has little to do with matrices, it being instead an immediate consequence of a fundamental abstract result; second, Oseledec's multiplicative ergodic theorem [26] (a different form of this was also proved by Millionshchikov [24]), which asserted in the matrix case that there moreover exists a random positive symmetric matrix $\Lambda$ for which

$$\lim_{n\to\infty} \frac{1}{n} \log \|(Z_n \Lambda^{-n})^{\pm 1}\| \to 0.$$

(This was not quite the original formulation; see [15].) In words, this says that there exists an "average" matrix $\Lambda = \Lambda(\omega)$ whose powers approximate the random product $Z_n(\omega)$, similar to the classical law of large numbers, except that $Z_n$ is written multiplicatively.

For a general group $G$, we need a substitute for the matrix norm. In the case of finitely generated groups, a good candidate is the word metric; see below. Ultimately, we will consider the situation of any $G$ which acts by isometry on a metric space. This includes the classical law of large numbers (translations are isometries of $\mathbf{R}$), Oseledec's theorem [invertible matrices are isometries of the space of symmetric positive definite matrices $\text{Pos}_N(\mathbf{R})$] and the action of a finitely generated group on itself (the word metric is left invariant), or, what amounts to the same thing, random walks on the underlying Cayley graph. *Thus, the setting of $G$ acting by isometry on a metric space $X$ is a very general and natural one for extensions of the law of large numbers.*

We now recall the concepts of Cayley graphs and word metrics. When a group $G$ is generated by a set $S$, one can consider the associated *Cayley graph* $X(G,S)$: the vertices are the elements of $G$, and two vertices $g$ and $h$ are adjacent if and only if they differ by an element of $S$ on the right, so $g = hs^{\pm 1}$. The action of $G$ on itself by left translation is an action of graph automorphisms of $X$ and hence an isometric action with respect to the associated graph distance, which is often called the *word metric*.

The most familiar examples of Cayley graphs are the standard graph of $\mathbf{Z}^N$ (a lattice) and the free group $F_N$ (a $2N$-regular tree). Suppose that the set $S$ is finite and the distribution of $g_i$ is uniform on $S$. Then notice that $Z_n = g_1 g_2 \cdots g_n$ is a simple random walk on the graph $X(G,S)$.

In this paper we prove a rather general law of large numbers for random walks on groups. We actually work in a setting more general than i.i.d., namely the *stationary* or *ergodic* setting. More precisely, let $(\Omega, \mu)$ be a standard Borel space with $\mu(\Omega) = 1$ and $L: \Omega \to \Omega$ an ergodic measure-preserving transformation. Let $G$ be a topological group (e.g., any group with the discrete topology) with its Borel $\sigma$-algebra. Assume that $g: \Omega \to G$ is a measurable map and let

$$Z_n(\omega) = g(\omega)g(L\omega) \cdots g(L^{n-1}\omega).$$

To fix the terminology, we will refer to $Z_n$ as an *ergodic cocycle* in the general case, and in the i.i.d. case, $Z_n$ is a *random walk*.

Let $X$ be a proper metric space (i.e., closed bounded sets are compact) and fix a basepoint $x_0 \in X$. A *horofunction* $h: X \to \mathbf{R}$ is a limit function

$$h(z) = \lim_{n \to \infty} d(x_n, z) - d(x_n, x_0)$$

for some sequence of points $x_n \to \infty$ in $X$, where the convergence is uniform on compact sets. Suppose $G$ acts by isometry on $X$ so that $\phi: G \to \mathrm{Isom}(X)$ is measurable (we will suppress $\phi$). For example, $X$ could be a Cayley graph of $G$ in the case where $G$ is finitely generated, or the space of positive definite symmetric matrices $\mathrm{Pos}_N$ in the case where $G$ is a group of real matrices. Assume that the cocycle is *integrable*, that is

$$\int_\Omega d(g(\omega)x_0, x_0) \, d\mu < \infty.$$

In the random walk case, this condition is that of *finiteness of the first moment*.

By subadditivity, which follows from the triangle inequality and the isometry property, Kingman's theorem implies that

$$A := \lim_{n \to \infty} \frac{1}{n} d(Z_n x_0, x_0)$$

exists almost surely and is independent of $\omega$, by ergodicity. We prove the following, general, noncommutative law of large numbers, or in a different terminology, multiplicative ergodic theorem:

THEOREM 1.1. *Let $X$ be a proper metric space and $Z_n$ an integrable ergodic cocycle taking values in $\mathrm{Isom}(X)$. Then, for almost every $\omega$, there is a horofunction $h = h_\omega$ depending measurably on $\omega$ such that*

$$\lim_{n \to \infty} -\frac{1}{n} h(Z_n x_0) = A,$$

*where $A := \lim_{n \to \infty} \frac{1}{n} d(Z_n x_0, x_0)$.*

This theorem makes a nontrivial statement about the behavior of simple random walk on *any* finitely generated nonamenable group (as it is known in this case that $A > 0$, see Theorem 1.3 below).

Here is a related result which applies to any group, in the sense that no properness is assumed:

THEOREM 1.2. *Assume that $Z_n$ is a ballistic random walk on a group. Then, for any $\varepsilon > 0$ and almost every trajectory, there is a time after which any finite collection of $\varepsilon$-shadows of the trajectory intersect.*

(See Section 5 for further explanation.) We will note at the end of Section 5 that, in fact, Theorem 1.1 implies Theorem 1.2 in the case where the metric space is proper. It is relevant to recall here a theorem of Guivarc'h [12]:

THEOREM 1.3. *Let $G$ be a locally compact group generated by compact set $V$ and denote by $\delta_V$ the corresponding word metric. Assume $\nu$ is a probability measure on $G$ whose support generates $G$ as a semigroup and denote by $Z_n$ the corresponding random walk. If $G$ is nonamenable and*

$$\int_G \delta_V(g)\, d\nu(g) < \infty,$$

*then there is a number $A > 0$ such that almost surely*

$$\lim_{n\to\infty} \frac{1}{n} \delta_V(Z_n) = A.$$

Guivarc'h wrote in [13] that this theorem does not treat the directional behavior of $Z_n$ and hence often gives very partial information. Our Theorems 1.1 and 1.2 do indeed provide information on the directional behavior of $Z_n$ and they do so even in the ergodic setting (whenever $A > 0$).

**2. Horofunctions.** Horofunctions have their origin in non-Euclidean geometry. More general definitions have been considered by Busemann and later by Gromov in [4], which was recalled in the Introduction. Horofunctions have been studied for spaces of nonpositive curvature, nonnegative curvature, as well as Gromov hyperbolic spaces; see [6]. More recent investigations on horofunctions include [28] and [21] for finite-dimensional Banach spaces and [30] for certain finitely generated groups.

A *geodesic ray* $\gamma$ is a map $\gamma : [0, \infty) \to X$ which is an isometry onto its image. Busemann associated to any geodesic ray $\gamma : [0, \infty) \to X$ with $\gamma(0) = x_0$, a horofunction $h_\gamma$, as follows:

$$h_\gamma(z) = \lim_{n\to\infty} d(\gamma(n), z) - n.$$

The limit exists by monotonicity coming from the triangle inequality. The convergence is moreover uniform if $X$ is proper, as can be seen from a $3\varepsilon$-proof using the compactness of closed balls.

More generally, let $X$ be a proper metric space (i.e., closed bounded sets are compact), and fix a basepoint $x_0 \in X$. Let

$$\Phi : X \to C(X)$$

be defined by $x \mapsto d(x, \cdot) - d(x, x_0)$ and where the topology on $C(X)$ is uniform convergence on compact sets. It can be checked that $\Phi$ is a continuous injection and we identify $X$ with its image. Let $H = \overline{\Phi(X)}$. It is easy to verify, since $X$ is proper and $|h(z)| \leq d(z, x_0)$, that $H$ is a compact metrizable space. By definition, the points in $H \setminus \Phi(X)$ are the horofunctions (based at $x_0$). The space $H$ is a compactification of $X$, since $X$ sits (via $\Phi$) inside it as an open dense subset.

The action of $\mathrm{Isom}(X)$ on $X$ extends continuously to an action by homeomorphisms to the whole of $H$ and is given by

$$g \cdot h(z) = h(g^{-1}z) - h(g^{-1}x_0),$$

as is straightforward to check.

A complete metric space $X$ is a CAT(0)-*space* or *nonpositively curved* if for any $x, y \in X$, there exists a point $z$ such that

$$d(x, y)^2 + 4d(z, w)^2 \leq 2d(x, w)^2 + 2d(y, w)^2$$

holds for every $w \in X$. This inequality is called the *semiparallelogram law*, motivated by the fact that in case of equality, it is the usual parallelogram law for Hilbert spaces. Apart from Euclidean spaces, other main examples are the classical hyperbolic spaces and $\mathrm{Pos}_N(\mathbf{R})$.

For simplicity, we assume in addition that $X$ is proper. There is a standard boundary $\partial X$ and compactification $X \cup \partial X$ associated to $X$ where $\partial X$ is the set of all geodesic rays from $x_0$. The topology is given by shadows of balls (cf. Section 5). For Euclidean or hyperbolic spaces, $\partial X$ is the sphere at infinity. This compactification turns out to be homeomorphic to $H$ above. We refer to [6] for more details on these topics.

**3. A law of large numbers for noncommuting random products.** Any group $G$ clearly acts faithfully by isometry on some metric space, since it, for example, acts on $l^2(G)$ (with respect to a Haar measure or the counting measure) or the cone of positive functions on $G$ equipped with Hilbert's metric. Many important groups, for example, every linear group, admit actions on nicer metric spaces than the ones just mentioned. Moreover, any finitely generated group acts faithfully on a proper metric space, since it acts by isometry on an associated Cayley graph. In view of this, we propose the following (already formulated in the Introduction) as an appropriate generalization of the classical strong law of large numbers to groups:

THEOREM 3.1. *Let $X$ be a proper metric space and $Z_n$ an integrable ergodic cocycle taking values in* $\mathrm{Isom}(X)$. *Then, for almost every $\omega$, there is a horofunction $h$ (depending on $\omega$) such that*

$$\lim_{n \to \infty} -\frac{1}{n} h(Z_n x_0) = A,$$

*where $A = \lim_{n \to \infty} d(Z_n x_0, x_0)/n$.*

When specialized to $G = GL_N(\mathbf{R})$ and $X = \mathrm{Pos}_N(\mathbf{R})$, this statement is in fact equivalent to Oseledec's theorem: $h$ corresponds to $\Lambda$ above. This is explained in [20]. It can be illustrated in the simplest case of $\mathbf{R}$, that is, the law of large numbers itself (or Birkhoff's ergodic theorem). Namely, the strong law of large numbers not only asserts that $|Z_n|/n \to A = |E(X_1)|$, but that actually $Z_n/n \to E(X_1)$. This shows the separateness of the issues of speed (the existence of $A$) and direction (the existence of $h$), in the most simple case.

QUESTIONS. What about random walks on graphs which are not homogeneous? What remains true if $X$ is no longer proper? Do limits along trajectories with respect to other horofunctions also exist? It seems plausible that to each of these questions, possibly with the exception of extensions to nonproper spaces, there will be counterexamples to the naïve extensions. What about Conjecture 8.1 in [20] where isometries are replaced by semi-contractions?

REMARK 3.2. The theorem also applies to the corresponding statement for Brownian motion $B_t$ (replace $Z_n x_0$ by $B_t$ in the statement) on Riemannian manifolds which have a cocompact group of isometries $\Gamma$, since one can then pass either to a random walk (Furstenberg, Lyons–Sullivan) or a ergodic cocycle on $\Gamma$ as is well known (see, e.g., Section 4 in [20]).

In order to explain which were the previously-known cases of the theorem, we first need to establish a proposition:

PROPOSITION 3.3. *Let $x_n$ be a sequence of points in $X$ and $A \geq 0$. Assume that there is a geodesic ray $\gamma$ such that $d(x_n, \gamma(An))/n \to 0$. Then*

$$\lim_{n \to \infty} -\frac{1}{n} h_\gamma(x_n) = A.$$

PROOF. For any horofunction $h$ it is true that $|h(x_n)| \leq d(x_n, x_0)$, from the triangle inequality. Since $d(x_n, x_0)/n \to A$, this implies that

$$\liminf_{n \to \infty} \frac{1}{n} h_\gamma(x_n) \geq -A.$$

On the other hand, note that for $t > An$, we have $d(\gamma(t), x_n) \leq t - An + d(\gamma(An), x_n)$. Hence,

$$h_\gamma(x_n) \leq -An + d(\gamma(An), x_n),$$

and the proposition follows on dividing by $n$ and taking the limit as $n \to \infty$. □

The notion of $\{x_n\}$ being of sublinear distance from a geodesic ray—*ray approximation*—was introduced and studied by Kaimanovich [14], who, in the case of symmetric spaces of nonpositive curvature (e.g., classical hyperbolic spaces and $\text{Pos}_N$), in fact characterized such sequences [15]. Hence Theorem 3.1 was known (in view of the proposition) in these cases. For general proper (also nonproper) CAT(0)-spaces, the theorem was established by Karlsson and Margulis in [19]. More precisely, an equivalent version of it was established (see [20]). There are some works on laws of large numbers in a stronger sense; see [16] concerning splittable solvable Lie groups. Note also that with the help of an idea of Delzant, Kaimanovich established ray approximation (and hence also Theorem 3.1 in this case) for Gromov hyperbolic spaces; see [18]. Some other papers establishing this type of law of large numbers include [7] and [25].

REMARK 3.4. In the case where $X$ is Gromov hyperbolic, it is known that for any two horofunctions $h_1$ and $h_2$ associated to sequences converging to $\xi$ in the Gromov boundary, there is a constant $C$ such that $|h_1(z) - h_2(z)| < C$ for all $z \in X$ (see [6]), and for a boundary point $\xi$, there may indeed be several such associated horofunctions. This shows that the $h$ in the theorem is not necessarily unique (it is unique only up to suitable equivalence).

**4. The proof of Theorem 1.1.** Assume the notation of the previous sections. In particular, let $X$ be a proper metric space on which $G$ acts by isometry and let $H$ be the compactification adding horofunctions to $X$.

Define for $g \in G$ and $h \in H$, $F(g, h) = -h(g^{-1}x_0)$ and note the following cocycle relation:

$$\begin{aligned} F(g_1, g_2 h) + F(g_2, h) &= -(g_2 \cdot h)(g_1^{-1} x_0) - h(g_2^{-1} x_0) \\ &= -h(g_2^{-1} g_1^{-1} x_0) + h(g_2^{-1} x_0) - h(g_2^{-1} x_0) \\ &= F(g_1 g_2, h). \end{aligned}$$

Note also that for any $g \in G$,

$$d(x_0, g x_0) = \max_{h \in H} F(g, h),$$

since $F(g, \Phi(g^{-1} x_0)) = -d(g^{-1} x_0, g^{-1} x_0) + d(g^{-1} x_0, x_0)$.

Let $Z_n(\omega)$ be an integrable cocycle taking values in $G$, as in the Introduction. We define the skew product $\overline{L} : \Omega \times H \to \Omega \times H$ via

$$\overline{L}(\omega, h) = (L\omega, g(\omega)^{-1} h).$$

Let $\overline{F}(\omega, h) = F(g(\omega)^{-1}, h)$. Since $|F(g(\omega)^{-1}, h)| \leq d(x_0, g(\omega) x_0)$ and because of the basic integrability assumption in the Introduction, $\overline{F}$ is in $L^1(\Omega, C(H))$. For detailed information about skew products that we will need here, see [1], Chapter 1.

Using the cocycle relation, we have

$$\overline{F_n}(\omega, h) := \sum_{i=0}^{n-1} \overline{F}(\overline{L}^i(\omega, h))$$
$$= F(g(\omega)^{-1}, h) + F(g(L\omega)^{-1}, g(\omega)^{-1} \cdot h) + \cdots$$
$$\quad + F(g(L^{n-1}\omega)^{-1}, g(L^{n-2}\omega)^{-1} \cdots g(\omega)^{-1} \cdot h)$$
$$= F(Z_n(\omega)^{-1}, h).$$

By the subadditive ergodic theorem,

$$A := \lim_{n \to \infty} \frac{1}{n} d(Z_n(\omega) x_0, x_0) = \inf_{n > 0} \frac{1}{n} \int_\Omega d(Z_n(\omega) x_0, x_0) \, d\mu(\omega)$$
$$= \inf_{n > 0} \frac{1}{n} \int_\Omega \max_{h \in H} \overline{F_n}(\omega, h) \, d\mu(\omega).$$

Consider the space of probability measures $\nu$ on $\Omega \times H$ which projects onto $\mu$ on $\Omega$, that is, $\nu(B \times H) = \mu(B)$ for any measurable set $B \subset \Omega$. The topology is the weak topology coming from the duality with $L^1(\Omega, C(H))$; see [1], page 27. For each $n$, choose a probability measure $\mu_n$ in this set such that

$$\frac{1}{n} \int_{\Omega \times H} \overline{F_n}(\omega, h) \, d\mu_n(\omega, h) \geq A.$$

For example, the measures defined by $\mu_\omega = \delta_{\Phi(Z_n(\omega))}$ in the terminology of [1], pages 22–25, would do.

Let

$$\eta_n = \frac{1}{n} \sum_{i=0}^{n-1} (\overline{L}^i)_* \mu_n$$

and let $\eta$ be a weak limit of these measures, which is possible due to the weak sequential compactness ([1], page 27).

The space of $\overline{L}$-invariant probability measures projecting onto $\mu$ and satisfying $\int \overline{F} \, d\nu \geq A$ is a compact, convex set. It is nonempty because $\eta$ belongs to this set. Indeed, it is clearly a probability measure, and the invariance is simple to check. Moreover, it is set up by construction so that $\int \overline{F} \, d\eta_n \geq A$ and, by definition of weak limits, this property passes to $\eta$ as well. The Krein–Milman theorem implies that the set has an extreme point $\eta_0$ which must be an ergodic measure for $\overline{L}$ by a standard argument.

The Birkhoff ergodic theorem implies that for $(\omega, h)$ in a set $P \subset \Omega \times H$ of $\eta_0$-measure 1, we have

$$\lim_{n \to \infty} \frac{1}{n} \sum_{i=0}^{n-1} \overline{F}(\overline{L}^i(\omega, h)) = \int_{\Omega \times H} \overline{F} \, d\eta_0 \geq A.$$

On the other hand, the left-hand side equals
$$\lim_{n\to\infty} -\frac{1}{n} h(Z_n(\omega)x_0) \leq A.$$
We therefore have equality everywhere. Since $\eta_0$ projects onto $\mu$, we have that for $\mu$-almost every $\omega$, there is a nonempty set of $h$ with the desired property. Finally, we will appeal to a measurable section theorem:

There is a Polish topology on $\Omega$ compatible with the standard Borel structure and such that the projection $f:\Omega \times H \mapsto \Omega$ maps open sets to Borel set, and the inverse image of each point in $\Omega$ is a closed subset. By regularity of $\eta_0$, we can find closed subsets of $P$ with arbitrarily large measure. These subsets are Polish spaces for the induced topology and the restriction of $f$ still satisfies the hypotheses of Theorem 3.4.1 in [2], which then gives a (partially-defined) cross section. Putting them together yields a measurable, a.e. defined cross section $\omega \mapsto (\omega, h_\omega)$ with $h_\omega$ having the desired property.

**5. Shadows of ballistic walks intersect.** Let $X$ be a metric space and fix a base point $x_0$. The $\varepsilon$-*shadow* of a point $y$ is
$$U_\varepsilon(y) := \{z : d(x_0, y) + d(y, z) \leq d(x_0, z) + \varepsilon d(x_0, y)\}.$$

In words, it is the set of points which almost lie on a geodesic ray passing through $x_0$ and $y$. More precisely, if $\gamma$ is a geodesic connecting $x_0$ and $y$, then $U_0(y)$ consists of all the points belonging to a geodesic extension of $\gamma$ beyond $y$. For CAT(0)-spaces, the sets $U_\varepsilon(z)$ constitute a basis of open sets for $X \cup \partial X$.

We will use a subadditive ergodic lemma from [19]. Note, however, that a lemma of similar type was proved by Pliss in [27] and that this has been very useful in smooth dynamics. Suppose $a(n,\omega)$ is a sequence of measurable functions satisfying
$$a(n+m,\omega) \leq a(n,\omega) + a(m, L^n\omega),$$
for every $n, m \geq 1$ and every $\omega$. Assume that
$$\int_\Omega |a(1,\omega)|\, d\mu(\omega) < \infty$$
and
$$A := \inf \frac{1}{n} \int_\Omega a(n,\omega)\, d\mu(\omega) > -\infty.$$
One can then prove (see [19], or [20] for an alternate proof):

LEMMA 5.1. *For almost every $\omega$, we have that for any $\varepsilon > 0$, there exists $K$ and infinitely many $n$ such that*
(1) $$a(n,\omega) - a(n-k, L^k\omega) \geq (A - \varepsilon)k$$
*for all $K \leq k \leq n$.*

Now let $G$ be a group. Assume that $G$ acts on a metric space $(X,d)$—for example, $X$ could be $G$ itself with an invariant metric $d$, such as a word metric. Let $g_1 : \Omega \to G \to \mathrm{Isom}(X)$ be a measurable map such that

$$\int_\Omega d(g_1(\omega)x_0, x_0)\, d\mu(\omega) < \infty.$$

Now, as usual, let $Z_n(\omega) = g_1(\omega)g_1(L\omega)\cdots g_1(L^n\omega)$ and assume that the cocycle or walk is *ballistic*, that is,

$$A := \inf_{n>0} \frac{1}{n} \int_\Omega d(Z_n(\omega)x_0, x_0)\, d\mu(\omega) > 0.$$

We can then prove the following:

THEOREM 5.2. *For almost every trajectory $\{Z_n(\omega)x_0\}$ and $\varepsilon > 0$, there is an integer $N = N(\omega, \varepsilon)$ such that for any $M > N$,*

$$\bigcap_{k=N}^M U_\varepsilon(Z_k x_0) \neq \varnothing.$$

PROOF. Let $\varepsilon > 0$ be given and let $\delta > 0$ be small so that $2\delta/(A+\delta) < \varepsilon$. Choose $N$ larger than $K$ in Lemma 5.1 applied to $a(n, \omega) := d(Z_n(\omega)x_0, x_0)$ with "$\varepsilon$" $= \delta$, and sufficiently large that $a(k, \omega) < (A+\delta)k$ for all $k \geq N$. By Lemma 5.1, there is an $n$ larger than $M$ such that

$$a(n, \omega) - a(n-k, L^k\omega) \geq (A-\delta)k$$

for all $N \leq k \leq n$. From this inequality, and in view of the invariance of $d$, we have

$$d(Z_n x_0, x_0) - d(Z_n x_0, Z_k x_0) \geq (A+\delta)k - 2\delta k$$
$$\geq \left(1 - \frac{2\delta}{A+\delta}\right)(A+\delta)k$$
$$\geq (1-\varepsilon)d(Z_k x_0, x_0).$$

Rearranging the terms yields

$$d(x_0, Z_k x_0) + d(Z_k x_0, Z_n x_0) \leq d(Z_n x_0, x_0) + \varepsilon d(x_0, Z_k x_0),$$

or, in other words, that $Z_n x_0 \in U_\varepsilon(Z_k x_0)$ for all $N \leq k \leq M$, and hence the intersection in the theorem is nonempty. □

REMARK 5.3. Note that the appearance of $\varepsilon$ is natural; $\varepsilon = 0$ would give very thin sets and the theorem would be false. The set $U_0(y)$ can be just a geodesic ray or even equal $\{y\}$ in general. One could use the shadows to contruct a boundary such that ballistic random walks would converge to points in this boundary. This is consistent with the fact that a simple, symmetric random walk on a finitely generated group where $A > 0$ has a nontrivial Poisson boundary, as proved by Varopoulos in [29].

REMARK 5.4. The same proof as above works if one replaces isometries by semicontractions, which, by definition, are self-maps $X \to X$ with Lipschitz constant 1.

Here is a corollary to illustrate the above phenomenon (a stronger result was proven in [19], with more elaborate geometric arguments, and in Theorem 1.1 above):

COROLLARY 5.5. *Assume that $X$ is a proper* CAT(0)-*space. Then, for almost every $\omega$, $Z_n(\omega)x_0$ converges to a point in $\partial X$.*

PROOF. Suppose there were two limit points $\gamma_1$ and $\gamma_2$ in $\partial X$. Then the $\varepsilon$-shadows, for some fixed $\varepsilon > 0$, of points in the subsequence approaching $\gamma_1$ versus the ones approaching $\gamma_2$, will eventually stop to intersect (in view of the fact that shadows generate the topology of $\partial X$). This contradicts the theorem. □

Apart from the symmetric spaces of nonpositive curvature, such as $\text{Pos}_N$, and spaces of pinched negative curvature, the corollary was first proved by Ballmann [3] in the case of random walks on $G$ acting cocompactly (or just satisfying the duality condition) and $X$ having rank 1. Note also that, in general, $A > 0$ is necessary, since the corollary is false for symmetric random walks on $\mathbf{Z}^N$ and $\mathbf{R}^N$.

Finally, we compare the theorems we obtain. Note that $z \in U_\varepsilon$ if and only if $\Phi(z)(y) \leq (\varepsilon - 1)d(x_0, y)$. Therefore, $V_\varepsilon(y)$, the closure of $\Phi(U_\varepsilon(y))$ in $H$, is the set of horofunctions $h$ such that $h(y) \leq (\varepsilon - 1)d(x_0, y)$. Theorem 1.1 says that for almost every $\omega$, the intersection over all sufficiently large $n$ of the sets $V_\varepsilon(Z_n(\omega))$ is not empty in $H$. This gives Theorem 5.2 (in the proper case) since the finite intersections have open interiors; if they are nonempty, they have to contain points from $\Phi(X)$.

**6. Comments on two special cases.** We wish to end this paper by recording a couple of related results which have not appeared in the literature.

6.1. *Coxeter groups.* Let $(W, S)$ be a finitely generated Coxeter group. These groups arise in several areas of mathematics. Moussong showed that there is an associated CAT(0)-space $X$ on which $W$ acts properly and cocompactly by isometries; see [6]. It follows that the number of orbit points grows at most exponentially with the radius. This is a case where Corollary 6.2 in [19] applies and we can state the following:

THEOREM 6.1. *Let $(W,S)$ be a finitely generated Coxeter group and $(X,d)$ its associated Moussong complex. Assume $\nu$ is a probability measure on $W$ such that*

$$\int_W d(gx_0, x_0)\, d\nu(g) < \infty$$

*for some (and hence any) $x_0 \in X$. Then the Poisson boundary of $(W,\nu)$ is either trivial or isomorphic to $\partial X$ with the induced hitting measure.*

6.2. *Mapping class groups.* The mapping class groups of surfaces play an important role in low-dimensional topology. They act by isometry on the associated Teichmüller space equipped with the Weil–Petersson metric. This metric is negatively curved but not complete. The main theorem in [19], however, assumes completeness of the space. But in view of recent works on the Weil–Petersson metric and a proof analysis, we can now formulate the following new theorem:

THEOREM 6.2. *Let $Z_n$ be an integrable ergodic cocycle taking values in the mapping class group $M_g$ and let $x_0$ be a point in $Teich_g$. Then almost every trajectory $Z_n x_0$ lies on sublinear distance from a Weil–Petersson geodesic ray $\gamma$, or equivalently,*

$$\lim_{n \to \infty} -\frac{1}{n} h_\gamma(Z_n x_0) = A,$$

*where $A = \lim_{n \to \infty} d(Z_n x_0, x_0)/n$ and $d$ is the Weil–Petersson metric.*

PROOF. The Weil–Petersson metric on Teichmüller space $Teich_g$ is known to have everywhere negative sectional curvature, but it is not complete. It is, however, geodesically convex and moreover, recent investigations (see [8]) on the metric space closure have shown that geodesics which meet the boundary (at finite distance) terminate and cannot be extended. In the proof of [19], it is only used that any two points can be joined by a geodesic segment in $X$. The limiting geodesic constructed there is the limit $\gamma(R)$ of Cauchy sequences $\gamma_i(R)$ of points on geodesics. The limit belongs to the completion. But now, in view of the above facts about Weil–Petersson geodesics, it is clear that all the points in $\gamma$ must actually lie inside the space $Teich_g$ (as opposed to merely in the completion). So the statement of the main theorem in [19] holds also in the current situation, which in turn is equivalent to the horofunction statement as proven in [20]. □

Note here that an identification of the Poisson boundary for random walks on the mapping class groups was obtained by Kaimanovich and Masur in [17]. Quite recently, Duchin [9] proved a multiplicative ergodic theorem (geodesic ray approximation) for the mapping class groups acting on the Teichmüller space with the Teichmüller metric.

**Acknowledgments.** The first author wishes to thank the Yale Mathematics Department for the invitation to spend the spring semester of 2005 there, during which part of this paper was written. Furthermore, he is grateful to Anna Erschler for inviting him to the stimulating Random Walks meeting in Lille, December 2004, and to Leonid Potyagailo for asking about the Weil–Petersson metric in the context of laws of large numbers.

DEPARTMENT OF MATHEMATICS
ROYAL INSTITUTE OF TECHNOLOGY (KTH)
100 44 STOCKHOLM
SWEDEN
E-MAIL: karlsson@aya.yale.edu

DEPARTMENT OF MATHEMATICS
HURLEY HALL
UNIVERSITY OF NOTRE DAME
NOTRE DAME, INDIANA 46556
USA
E-MAIL: fledrapp@nd.edu